\documentclass[12pt,leqno]{amsart}
\usepackage{amsmath}
\usepackage{graphicx,color}
\usepackage{verbatim}
\newtheorem{theorem}{Theorem}[section]

\newtheorem{lemma}[theorem]{Lemma}
\newtheorem{proposition}[theorem]{Proposition}
\theoremstyle{definition}
\newtheorem{definition}[theorem]{Definition}
\theoremstyle{remark}
\newtheorem{remark}[theorem]{Remark}

\numberwithin{equation}{section}

\def\R{\mathbb{R}} 
\def\C{\mathbb{C}} 

\def\rt{\rightarrow}

\begin{document}

\sloppy

\title{Totally Geodesic discs in strongly convex domains}

\author{Herv{\'e} Gaussier and Harish Seshadri}

\maketitle

\section{Introduction}
The Kobayashi metric and its infinitesimal version, introduced by S.Kobayashi \cite{kob1}, carry geometric properties of complex manifolds.
In case the integrated version of the infinitesimal metric, or equivalently the Kobayashi metric (see \cite{roy}) define a distance, called the Kobayashi distance, the associated metric space inherits dynamical and geometric properties fitted to the study of holomorphic function spaces and the associated metric space is named Kobayashi hyperbolic. One may refer to \cite{kob2} for a general presentation of Kobayashi hyperbolic spaces. Strictly pseudoconvex domains in the complex Euclidean space $\mathbb C^n$ or more restrictively strongly convex domains, endowed with their Kobayashi distance, are classical important examples of Kobayashi hyperbolic spaces. Let $\Omega$ be a bounded $C^3$ strongly convex domain in $\C^n$, i.e.,  for any two points $z,z' \in \overline{\Omega}$ the open line segment $(z,z')$ is contained in $\Omega$. In this paper we investigate some geometric aspects of the metric space $(\Omega, d_\Omega^K)$, where $d^K_\Omega$ denotes the Kobayashi metric of $\Omega$.
We recall that a map $f : \mathbb C^{n_1} \to \mathbb C^{n_2}$ is anti-holomorphic if $df \circ J_{1} = -J_{2}\circ df$ where $J_{i}$ denotes the standard complex structure on the Euclidean space $\mathbb C^{n_i},\ i=1,2$. Let $\Delta$ denote the unit disc in $\C$. Our main result is then the following:
\begin{theorem}\label{mai1}
Let $f: (\Delta, d^K_\Delta) \rt (\Omega, d^K_\Omega)$ be an isometry, namely
$$d^K_\Omega(f(\zeta),f(\eta)) = d^K_\Delta (\zeta,\eta) \ \ \ \ \ \forall  \ \zeta,\eta \in \Delta .$$

Then $f$ is either holomorphic or anti-holomorphic.
\end{theorem}

As a corollary we have
\begin{theorem}\label{mai2}
Let $n_1, n_2$ be positive integers and let $\Omega_i \subset \C^{n_i}$, $i=1,2$, be bounded $C^3$ strongly convex domains. If $\phi: (\Omega_1, d^K_{\Omega_1}) \rt (\Omega_2, d^K_{\Omega_2})$ is an isometry then $\phi$ is either holomorphic or anti-holomorphic.
\end{theorem}

We note that no assumptions are made about the smoothness of the isometry.  \vspace{2mm}

In ~\cite{kv1} an analogue of Theorem \ref{mai2} was proved under the hypotheses that the domains are equidimensional and strongly pseudoconvex but with the stronger assumption that the isometry is $C^1$ and has a $C^1$ extension to the boundary of $\Omega_1$. In
~\cite{kv2} an analogue of the Wong-Rosay theorem about noncompact automorphism groups was proved for equidimensional strongly convex domains. This was extended to strongly pseudoconvex domains in ~\cite{kk}. To the best of our knowledge, the question of whether an isometry between strongly pseudoconvex domains (even in the equidimensional case)  is holomorphic or anti-holomorphic is still open.

The proof of Theorem \ref{mai1} proceeds as follows. We first observe that any isometric map $\gamma: I \rt (\Omega, d^K_\Omega)$, where $I \subset \R$ is an interval, is a real geodesic i.e., a Kobayashi length minimizing $C^1$ curve. In fact we prove that any such map in $\Omega$ is contained in a complex geodesic (in the sense of Lempert).  Let $f: \Delta \rt \Omega$ be a $C^1$ isometry. Choose two real geodesics $\sigma$ and $\gamma$ in $\Delta$ which approach the same point $w \in \partial \Delta$. We can reparametrize these geodesics to get $\sigma_1$ and $\gamma_1$ which are now smoothly defined on $[0,1]$ with $\sigma_1 (1)= \gamma_1 (1) =w$. Since  $f \circ \sigma$ and  $f \circ \gamma$ are isometric maps of intervals into $\Omega$ they are smooth. Moreover the corresponding reparametrizations $f \circ \sigma_1$ and $f \circ \gamma_1$ also extend smoothly to $[0,1]$. We then prove the key fact that  $(f \circ \sigma_1)'(1) = (f \circ \gamma_1)'(1)$.

Let $\phi$ and $\psi$ be two complex geodesics such that after composing with an automorphism of $\Delta$ or with the conjugate of an automorphism of $\Delta$ we get
$\phi \circ \sigma = f \circ \sigma$ and $\phi \circ \gamma = f \circ \gamma$. One can then see that $(d \phi)(w)= (d \psi)(w)$. On the other hand, we prove that if two complex geodesics agree up to first order at a boundary point then their images coincide. We point out that a similar result (about the uniqueness of complex geodesics with prescribed boundary data) was proved in \cite{chl} under the stronger assumption that $\Omega$ is of class $C^{14}$. We finally proved that $f(\Delta)=\phi(\Delta)$. It remains to prove that $df$ commutes (or anti commutes) with the standard structures on $\Delta$ and $\mathbb C^n$. That completes the proof of Theorem \ref{mai1}.
\vspace{2mm}

Theorem \ref{mai2} is a direct corollary of Theorem \ref{mai1}. The proofs of Theorem \ref{mai1} and Theorem \ref{mai2} are given in Section 4.
\vspace{2mm}

Finally we make a few remarks about our terminology. A smooth embedding between Riemannian manifolds
$f:(M,g) \rt (N,h)$ is said to be {\it totally geodesic} if $f^{\ast}(h)=g$ and the second fundamental form
of the image $f(M)$ vanishes. It can be checked that a smooth embedding is totally geodesic if and only if
it is a metric space isometry from $(M,d_g)$ to $(N,d_h)$ where $d_g$ and $d_h$ are the distance functions induced from $g$ and $h$ respectively. Hence we use the terms ``totally geodesic map" and ``isometry" interchangeably.

\section{Preliminaries}
 Given a bounded domain $D\subset \mathbb C^n$ we denote by $K_D$ the Kobayashi infinitesimal pseudometric on $D \times \mathbb C^n$, by $d^K_D$ the Kobayashi distance on $D$. In case $D=\Delta$, where $\Delta$ is the unit disc in $\mathbb C$, then $d^K_\Delta$ is the Poincar\'e distance on $\Delta$.

In this Section we collect some basic facts about the geometry of the Poincar\'e disc $(\Delta,d^K_\Delta)$ and about the behaviour of complex geodesics in strongly convex domains in $\mathbb C^n$.

\begin{definition}\label{rg}
A {\it geodesic} (or {\it real geodesic}) in $\Omega$ is a smooth curve $\gamma : I \rt \Omega$ such that
 $$l(\gamma \vert_{[t_1,t_2]}) := \int_{t_1}^{t_2} K_\Omega (\gamma(t), \gamma'(t))dt=d_\Omega^K(\gamma(t_1), \gamma(t_2))$$
 for all $t_1 , t_2 \in I$, where $I \subset \R$ is an interval.
\end{definition}


\subsection{Geometry of the Poincar\'e disc}
For $p \in \Delta$, let
$\delta(p)= dist(p, \ \partial \Delta).$
Note that
\begin{equation}\label{klm}
 d^K_\Delta(p,0)= -log\delta (p)
 \end{equation}
where $0 \in \Delta$ is the origin.

\begin{lemma}\label{dd}
Given $\epsilon >0$ there exists $C=C(\epsilon)>0$ such that the following holds:

Let $p, q \in \Delta$  satisfy $d^k_\Delta(p,q) \le \epsilon$. Then
$$C^{-1} \delta(p) \le \delta (q) \le C \delta (p).$$
Moreover $C \rt 1$ as $\epsilon \rt 0$.
\end{lemma}

\noindent {\bf Proof:} This follows from (\ref{klm}) since
$$d^K_\Delta(p,0) - \epsilon \ \le \ d^K_\Delta(q,0) \ \le  \ d^K_\Delta(p,0) + \epsilon.$$ \hfill $\square$

The following fact is standard.
\begin{lemma}\label{qq}
Let $\gamma, \sigma:[0,\infty) \rt \infty$ be unit-speed geodesics so that $\lim_{t \rt \infty} \gamma (t)=
\lim_{t \rt \infty} \sigma (t)=z \in \partial \Delta$.

Suppose that $\sigma(0)$ and $\gamma(0)$ lie on the same horocycle passing through $z \in \partial \Delta$. Then there is a constant $D$ such that
$$ d^K_\Delta (\gamma(t), \sigma(t)) \le D e^{-t} .$$
\end{lemma}

\subsection{Complex geodesics and holomorphic retracts in strongly convex domains}
Let $\phi : \Delta \to \Omega$ be a holomorphic disc.

\begin{definition}\label{cg}
 $(a)$ We call $\phi$ a {\sl  complex geodesic} if $\phi$ is an isometry for
the Kobayashi distances on $\Delta$ and $\Omega$. \vspace{2mm}

$(b)$ We call $\phi$ {\sl extremal with respect to $p,q \in \Omega$} if $\phi(0)=p$, $\phi(\zeta)=q$ for some $0 < \zeta < 1$ and $d^K_\Omega(p,q)=\log\{(1+\zeta)/(1-\zeta)\}$. \vspace{2mm}

$(c)$ We call $\phi$ {\sl extremal with respect to $(p,v) \in \Omega \times \mathbb C^n$} if $f(0)=p$, $df(0)v=\lambda v$, $\lambda > 0$, and if for every $\psi:\Delta \to \Omega$ such that $\psi(0)=p$, $d\psi(0)v=\mu v$ with $\mu > 0$, we have $\mu \leq \lambda$. \vspace{2mm}

$(d)$ A subset $S$ of a domain $D \subset \mathbb C^n$ is called a holomorphic
retract if there is a holomorphic mapping $r: D \to D$ such that $r(D) \subset S$ and $r(z)=z$ for $z \in S$.
\end{definition}

The following result due to L.Lempert is fundamental to this paper. Parts $(i)-(iv)$ are contained in \cite{lem1}, point $(v)$ is the content of Theorem 2 in \cite{lem2}. We point out that Theorem 2 in \cite{lem2} was stated for smooth $C^\infty$ domains but that the proof goes through for $C^3$ domains.

\begin{theorem}\label{lem-thm}
Let $\Omega$ be a bounded strongly convex domain in $\mathbb C^n$, with $\partial \Omega$ of class $C^3$.

(i) A map $\phi: \Delta \rt \Omega$ is a complex geodesic if and only if it is extremal with respect to any $(p,q) \in \phi(\Delta) \times \phi(\Delta)$ or with respect to any $(p,v) \in \phi(\Delta) \times \mathbb C^n$ (after composition with an automorphism of $\Delta$).

(ii) Given two points $p,q \in \Omega$ there is a unique complex geodesic $\phi$ whose image contains  $p$ and $q$.

(iii) Given a point $p \in \Omega$ and a 2-dimensional $J$-invariant subspace $V$ of $T_p \Omega$ (equivalently, a
complex tangent vector at $p$), there is a unique complex geodesic $\phi$ passing through $p$ and satisfying
$T_p (\phi (\Delta))=V$.


(iv) The map $\phi$ is proper and $\phi$ extends as a $C^1$ map up to $\overline{\Delta}$. Also $\phi(\overline \Delta)$ intersects $\partial \Omega$ transversally, namely $T^\mathbb C_{\phi(e^{i\theta})}(\partial \Omega) \cap T_{\phi(e^{i\theta})}(\phi(\Delta)) = \{0\}$.

(v) The one-dimensional holomorphic retracts in a strongly convex bounded domain are precisely the extremal discs.
\end{theorem}

Here $T^\mathbb C_{\phi(e^{i\theta})}(\partial \Omega)$ denotes the complex tangent space to $\partial \Omega$ at point $\phi(e^{i\theta})$.

\begin{lemma}\label{kol}
Let $\Omega$ be a $C^3$ strongly convex domain in $\C^n$ and let $\phi, \ \psi \ : \ \overline \Delta \rt  \overline \Omega$ be two complex geodesics. If there is a point $w \in \overline \Delta$ such that $\phi(w) = \psi(w)$ and $d \phi_w =d \psi_w$ then $\phi=\psi$.
\end{lemma}

\noindent{\bf Proof:} The result is immediate if $w \in \Delta$ since a complex geodesic passing through a point in a given direction is unique according to \cite{lem1}. Assume now that $w=1$, $\phi(1)=\psi(1)=1$, $d\phi_1= d \psi_1$ and $\phi \not\equiv \psi$ (notice that $\phi$ and $\psi$ are $C^1$ maps on $\overline{\Delta}$ by Theorem \ref{lem-thm}, point $(iv)$).
We keep the same notations as in the proof of Proposition 8 in \cite{lem1}. In particular for two elements $z=(z_1,\dots,z_n),\ w=(w_1,\dots,w_n) \in \mathbb C^n$ we set $\langle z,w\rangle := \sum_{j=1}^n z_jw_j$. If $z \in \partial \Omega$ denote by $\nu(z)$ the (outward) normal vector to $\partial \Omega$ at $z$.
According to \cite{lem1} there is a positive function $p$, continuous on $\partial \Delta$ such that the map
$\zeta \in \partial \Delta \mapsto \zeta p(\zeta)\overline{\nu(\phi(\zeta))}$ extends to a map $\tilde \phi$, continuous on $\overline{\Delta}$, holomorphic on $\Delta$ (see page 434 of \cite{lem1}).

Since $\Omega$ is strongly convex there is a constant $C>0$ such that :
$$
Re\langle \phi(\zeta)-\psi(\zeta),\overline{\nu(\phi(\zeta))}\rangle \geq C
$$
on a subset of positive measure in $\partial \Delta$.

Since $-(\zeta-1)^2/\zeta >0$ for $\zeta \in \partial \Delta \backslash \{1\}$ then, changing $C$ if necessary~:
$$
Re\left\langle \frac{\phi(\zeta)-\psi(\zeta)}{(\zeta-1)^2},\tilde{\phi}(\zeta) \right\rangle \geq C
$$
on a subset of positive measure in $\partial \Delta$.

Hence we have :

$$
Re\left(\frac{1}{2i\pi}\int_{\partial \Delta}\left\langle \frac{\phi(\zeta)-\psi(\zeta)}{(\zeta-1)^2},\tilde{\phi}(\zeta)\right\rangle \frac{d\zeta}{\zeta} \right)\geq C.
$$
However :
$$
\frac{1}{2i\pi}\int_{\partial \Delta} \left\langle \frac{\phi(\zeta)-\psi(\zeta)}{(\zeta-1)^2},\tilde{\phi}(\zeta)\right\rangle \frac{d\zeta}{\zeta} = \langle \phi(0)-\psi(0),\tilde{\phi}(0)\rangle,
$$
since the maps $\tilde{\phi}$ and $\displaystyle \zeta \in \overline{\Delta} \mapsto \frac{\phi(\zeta)-\psi(\zeta)}{(\zeta-1)^2}$ are holomorphic on $\Delta$ and continuous on $\overline \Delta$.
Hence :
\begin{equation}\label{geq}
Re\langle \phi(0)-\psi(0),\tilde{\phi}(0)\rangle \geq C.
\end{equation}

Following the proof of Proposition 2 in \cite{lem1}, for $\eta \in \Delta$, let $a_\eta$ be the automorphism of $\Delta$ defined by $a_\eta(\zeta) = (\zeta + \eta)/(1+\overline{\eta}\zeta)$. Since the index of the function $\zeta \in \partial \Delta \mapsto \zeta/a_{\eta}(\zeta)$ is not zero on $\partial \Delta$ we may choose a holomorphic function $q_\eta$ on $\Delta$ such that $Im(q_{\eta}(\zeta))=Im(\log(\zeta/a_\eta(\zeta)))$ for $\zeta \in \partial \Delta$. Note that we may fix the value $Re q_\eta(0)=0$ for every $\eta$.


We may apply the inequality~(\ref{geq}) to $\phi \circ a_\eta$ and $\psi \circ a_\eta$. This gives for every $\eta \in \Delta$ :

$$
Re \langle \phi(\eta)-\psi(\eta),\tilde{\phi_\eta}(0)\rangle \geq C.
$$
Here, according to the proof of Proposition 2 in~\cite{lem1}, $\tilde{\phi_\eta}(\zeta)=\exp\left(q_\eta(\zeta)\right) \tilde{\phi}(a_\eta(\zeta))$ for every $\zeta \in \Delta$. In particular $\tilde{\phi_\eta}(0)=\exp(i Im(q_\eta(0)))\tilde{\phi}(\eta)$.

We finally proved for every $\zeta \in \Delta$ :

$$
Re \langle \phi(\eta)-\psi(\eta),\exp(i Im(q_\eta(0)))\tilde{\phi}(\eta)\rangle \geq C.
$$
This is a contradiction for $\eta \to 1$ since $\tilde \phi$ is continuous on $\overline \Delta$ and $\phi(1)=\psi(1)$. \qed

\section{Totally geodesic discs}
We begin by noting that a real geodesic in $\Omega$ is an isometry from $I$ to $\Omega$ where $I$ carries the usual Euclidean distance. We first prove the partial converse that isometries from an interval $I \subset \R$ to $(\Omega, d^K_\Omega)$ are absolutely continuous and their lengths realize Kobayashi distance.
\begin{lemma}\label{ae}
Let $I \subset \R$ be an interval and $\alpha: I \rt \Omega$ an isometry. Then $\alpha$ is locally Lipschitz.
In particular $\alpha$ is absolutely continuous and
$$\int_{t_1}^{t_2} K_\Omega (\gamma(t), \gamma'(t))dt=d_\Omega^K(\gamma(t_1), \gamma(t_2)) \ \ \ \ \ \forall t_1,t_2 \in I. $$
\end{lemma}
{\bf Proof:} Without loss of generality assume that $I=[0,l]$ for some $ l>0$.
Let $p = \alpha(0)$ and consider the ball $B=B(p, 2l)$ with center $p$ and radius $2l$ in
the Kobayashi metric $d^K_\Omega$. By continuity of the infinitesimal Kobayashi metric there exists
$C>0$ such that
$$C \Vert v \Vert \le K_\Omega (q,v) \ \ \ \ \ \forall \ q \in B, \  v  \in \C^n.$$
Let $t_1, t_2 \in [0,l]$. Integrating the above estimate along the geodesic $\alpha \vert_{[t_1,t_2]}$ we get
$$C \Vert \alpha(t_2) - \alpha(t_1) \Vert \le d^K_\Omega(\alpha(t_2), \alpha(t_1))= \vert t_2 - t_1 \vert.$$
This proves that $\alpha$ is Lipschitz.

To see the second part, we note the following fact which is the content of Theorem 1.2 in \cite{ven}. If $\gamma : [0,l] \rt \Omega$ is an absolutely continuous curve then
\begin{equation}
 \int_0^l K_\Omega (\gamma(t), \gamma'(t)) dt = \sup_P \sum_{i=1}^{k-1}d^K_\Omega (\gamma(t_i), \gamma (t_{i+1})
\end{equation}
where the sum is over all partitions $P= \{ t_1=0,...,t_k=l\}$ of $[0,l]$.
If $\gamma$ is an isometry then
$$\sum_{i=1}^{k-1}d^K_\Omega (\gamma(t_i), \gamma (t_{i+1}))=d(\gamma(0), \gamma(l))$$
for any partition $P= \{ t_1=0,...,t_k=l\}$ of $[0,l]$ and the proof is complete. \hfill $\square$ \vspace{2mm}


\begin{definition}
We say that a map $f: \Delta \rt \Omega$ is a {\it totally geodesic disc} if $f$ is an isometry  for the Kobayashi distance :
$$
d^K_\Delta(x,y) = d^K_\Omega(f(x),f(y))
$$
for any two points $x,y \in \Delta$.
\end{definition}

 It can be checked that a  totally geodesic disc  $f$ is a proper map and extends to a map of class at least $C^{1/2}$ up to $\partial \Delta$. We will not prove these facts since we will not use them.

\begin{lemma}\label{geod-lem}

Let $I \subset \R$ be an interval and let $\alpha: I \rt \Omega$ be an isometry. Then the image of  $\alpha$ is contained in a complex geodesic i.e. there exists a complex geodesic $g : \Delta \to \Omega$ such that $\alpha (I) \subset g( \Delta )$.

In particular, every isometry $\alpha: I \rt \Omega$ is $C^1$ and there is a unique real geodesic between any two points in $\Omega$.

\end{lemma}

\noindent {\bf Proof:} Without loss of generality assume that
$I=[0,t_0]$ for some $t_0 >0$ and that $\alpha$ is differentiable at $t=0$. Let $\phi: \Delta \rt \Omega$ be a complex geodesic joining the points $p:=\alpha(0)$
and $q:=\alpha(t_0)$ and let $\gamma: [0,t_0] \rt \Omega$ be the real geodesic connecting $p$ and $q$ i.e. $\gamma(0)=p$, $\gamma(t_0)=q$, and contained in $\phi(\Delta)$. Let $\pi: \Omega \rt  \ \phi(\Delta)$ be the Lempert retract corresponding to $\phi$ (see Theorem \ref{lem-thm}, point $(v)$). We first note that
\begin{equation}\label{aaa}
 \pi \circ \alpha = \gamma.
\end{equation}
and
\begin{equation}\label{sss}
K_\Omega( \pi \circ \alpha (t), d \pi (\alpha'(t)))= K_\Omega(  \alpha (t), \alpha'(t))
\end{equation}
for almost all $t \in [0,t_0]$.

This is because $l( \pi \circ \alpha) = \int_{0}^{t_0} K_\Omega( \pi \circ \alpha (t), d \pi (\alpha'(t))) dt
 \le \int_{0}^{t_0} K_\Omega(  \alpha (t), \alpha'(t)) dt = l(\alpha)$. Here we have used the decreasing property of the Kobayashi norm under holomorphic mappings. Next we note that length minimizing curves are unique
 for the Poincar\'e metric on $\Delta$ and hence unique in $\phi(\Delta)$. Since $\pi \circ \alpha$ joins $p$ and $q$ and $l(\pi \circ \alpha) =l(\alpha)$ we get $ \pi \circ \alpha = \gamma$ and
$K_\Omega( \pi \circ \alpha (t), d \pi (\alpha'(t)))= K_\Omega(  \alpha (t), \alpha'(t))$ almost everywhere on $[0,t_0]$.

Next we claim that $\alpha'(0)=\gamma'(0)$. Let $N= {\rm Ker} (d \pi_p) \subset T_p \Omega$. Consider the function
$ f: N \rt [0, \infty)$ defined by
$$ f(n) = K_\Omega(p, v-n)$$
where $v =\alpha'(0)$. Since $\Omega$ is a strongly convex domain the Kobayashi indicatrix $I_p(\Omega):=\{x \in T_p\Omega /\ K_\Omega(p,x) < 1\}$ is strongly convex for every $p \in \Omega$ (see for instance \cite{pat}). The closure $\{x \in T_p\Omega /\ K_\Omega(p,v) \leq 1\}$ of $I_p(\Omega)$ is also strongly convex. It follows now from the homogeneity property $K_\Omega(p,cx) = |c|K_\Omega(p,x)$ (for every $x \in T_p\Omega,\ c \in \mathbb R$) that the set $\{x \in T_p\Omega /\ K_\Omega(p,x) \leq c\}$ is strongly convex for every $c>0$. Since $v+N$ is an affine subspace of $T_p\Omega=\mathbb C^n$ not containing the origin,
the function $f$ above attains its infimum $\inf f$ at exactly one point $n_0$.
Write $ v =n_0 + h_0$ and note that $K_\Omega(p,h_0) \le K_\Omega(p, v)$ by definition of $h_0$. On the other hand $d \pi(h_0)= d \pi (v)$ and $K_\Omega(p, h_0) \ge K_\Omega(p, d \pi(h_0))= K_\Omega(p, d \pi (v))= K_\Omega(p, v)$. Hence $K_\Omega(p,h_0) = K_\Omega(p, v)$ and $n_0=0$ by the uniqueness of the minimum of $f$.  If we let $n_1 = v -\gamma'(t_0)$ then $n_1 \in N$ by (\ref{sss}).  Moreover $K_\Omega(p, v-n_1)= K_\Omega(p, \gamma'(t_0))=K_\Omega(p, v)$ by (\ref{aaa}). Again by the uniqueness of minima of $f$, $n_1=0$ i.e.
$\gamma'(0)=v$.

Choose  any $t \in (0,t_0)$ such that $\alpha$ is differentiable at $t$ and consider the geodesic segment $\alpha \vert_{[0,t]}$. Let
$q_t= \alpha(t)$ and $\phi_t : \Delta \rt \Omega$ the complex geodesic passing through $p=\alpha(0)$ and $q_t$. Let $\gamma_t$ the
corresponding real geodesic connecting $p$ and $q_t$ which lies on the image of $\phi_t$.  The argument above applied to this new configuration gives  $\gamma_t'(0)=v$. The holomorphicity  of $\phi$ and $\phi_t$ imply that the tangent
spaces $T_p \phi(\Delta)=T_p \phi_t (\Delta)$. By Lemma~\ref{kol} we have  $\phi(\Delta)=\phi_t(\Delta)$. In particular $\alpha(t) \in \phi(\Delta)$. Since the set of points where $\alpha$ is differentiable has full measure and
$\phi(\Delta)$ is closed in $\Omega$, this completes the proof of the statement $\alpha(I) \subset \phi(\Delta)$.

By (\ref{aaa}) and the uniqueness of length minimizing curves in $\phi(\Delta)$ it follows that $\alpha=\gamma$ and the other statements of Lemma \ref{geod-lem} follow as well.  \hfill $\square$

\begin{remark}\label{geod-rem}
According to lemma \ref{geod-lem} let $\alpha$ be a real geodesic in $\Delta$ and let $f$ be a totally geodesic map in $\Omega$. Since $f \circ \alpha$ is a real geodesic in $\Omega$ there is a unique complex geodesic $g$ in $\Omega$ such that $f \circ \alpha((-\infty,\infty))$ is a smooth curve in $g(\Delta)$. Since $g$ is an embedding and an isometry for the Kobayashi metric, there is a unique real geodesic $\tilde \alpha$ in $\Delta$ such that $g \circ \tilde \alpha = f \circ \alpha$. Finally, after composing $g$ with an automorphism of $\Delta$ or with the conjugate of an automorphism of $\Delta$ , denoted by $\mu$, we may assume that $g \circ \mu \circ \alpha = f \circ \alpha$. We point out that the map $g \circ \mu$ is either holomorphic or anti-holomorphic.
\end{remark}
\vspace{3mm}

Given a real geodesic $\alpha:[0,\infty) \rt \Omega$ reparametrize $\alpha$ to get $\alpha_1: [0,1) \rt \Omega$, where
$$\alpha_1(u)= \alpha(-log(1-u)).$$
Then we have :

\begin{lemma}\label{lem-repara}
 $\alpha_1:[0,1) \rt \Omega$ extends $C^1$-smoothly to $[0,1]$ and it meets $\partial \Delta$ transversally (at $\alpha _1(1))$.
\end{lemma}

\noindent{\bf Proof:} According to Lemma \ref{geod-lem} consider the complex geodesic $g:\Delta \to \Omega$ such that $\alpha_1([0,1)) \subset g(\Delta)$. It follows from \cite{lem1} that $g$ extends to $\overline\Delta$ as a map of class $C^1$. Keeping the notations of Remark \ref{geod-rem} we may assume that $\alpha_1([0,1))=(g \circ \mu)([t_0,1))$ for some $-1 < t_0 < 1$. Hence $\alpha_1$ extends $C^1$-smoothly to $[0,1]$. The transversality of the intersection $\alpha_1([0,1])$ and $\partial \Omega$ is now a direct consequence of the estimates  of the Kobayashi infinitessimal metric on $\Omega$ (see \cite{gra}). \qed

\vskip 0,2cm
The following lemma is crucial for the results of this paper:

\begin{lemma}\label{lem-lambda}
Let $\gamma, \sigma : [0,\infty) \rightarrow \Delta$ be two  geodesics parametrized with respect to arc-length so that

\noindent (i) \ $\lim_{t \rightarrow \infty} \gamma(t) =\lim_{t \rightarrow \infty} \sigma (t)=1 \in \partial \Delta$

\noindent (ii) \ $\sigma(0)$ and $\gamma(0)$ lie on the same horocycle passing through $1 \in \partial \Delta$.

\noindent Let $\gamma_1, \ \sigma_1: [0,1] \rt \Omega$ be the corresponding reparametrizations. If $f: \Delta \rt \Omega$ is a totally geodesic disc then
$$ (f \circ \gamma_1)'(1)=(f \circ \sigma_1)'(1).$$

\end{lemma}


\noindent{\bf Proof:}
By Lemmas~\ref{dd} and \ref{qq}  there exists  constants $C_1, \ C_2 >0$ such that
$$ C_1^{-1} e^{-t} \ \le \ \delta( \gamma(t)) \ \le \ C_1 e^{-t},$$
and

 $$ d^K_\Delta (\gamma(t), \sigma(t)) \le C_2 e^{-t} .$$

Pick a sequence $t_\nu \rightarrow \infty$ and let
$x^\nu:=\gamma(t_\nu), \ y^\nu:=\sigma(t_\nu)$.
Then $\lim_{\nu \to \infty}x^\nu = \lim_{\nu \to \infty} y^\nu = 1$ and
$$\displaystyle \lim_{\nu \rt \infty} \frac{d^K_\Delta(x^\nu,y^\nu)}{\sqrt{\delta(x^\nu)}}=0.$$



According to Lemma \ref{geod-lem}, since $f \circ \gamma (-\infty,\infty)$ is a real geodesic in $\Omega$, there is a unique complex geodesic $g$ contained in $\Omega$ such that $f \circ \gamma((-\infty,\infty))\subset g(\Delta)$. According to Remark \ref{geod-rem} there is an automorphism $\mu$ of $\Delta$ (or the conjugate of an automorphism of $\Delta$) such that $(g \circ \mu)(\gamma (0))=f(\gamma(0))$.
Moreover it follows from the Hopf Lemma applied to the complex (or anti-complex) geodesic $g \circ \mu$ that the Euclidean distances $\delta(g \circ \mu(x^\nu))$ and $\delta(x^\nu)$ are equivalent. We finally have :
\begin{equation}\label{equiv-eq}
\displaystyle \lim_{\nu \rt \infty} \frac{d^K_\Omega(f(x^\nu),f(y^\nu))}{\sqrt{\delta(f(x^\nu))}} =0.
\end{equation}

Since $f$ extends up to $\partial \Delta$ we know that $f(1)=y^\infty \in \partial \Omega$. \vspace{2mm}

\noindent {\bf Claim:} $(f \circ \gamma_1)'(1) = \lambda (f \circ \sigma_1)'(1)$ for some $\lambda \in (0,\infty)$. \vspace{2mm}

This claim will follow from the following two lemmas. \vspace{2mm}

\begin{lemma}\label{qqq}
For any geodesic $\gamma: [0,\infty) \rt \Delta$ we have $\lim_{u \rt 1}(f \circ \gamma_1)'(1) \neq 0$.
\end{lemma}
\noindent {\bf Proof:} We have
$$  \gamma_1'(u) = (1-u)^{-1} \gamma'(t)$$
where $t=-log(1-u)$, $u \in [0,1)$.

Now
\begin{equation} \label{pop}
(f \circ \gamma)'(t)= df(\gamma(t))(\gamma'(t))= (1-u) df(\gamma_1(u))\gamma_1'(u)=(1-u) (f \circ \gamma_1)'(u).
\end{equation}

Since
$$K_\Omega (f \circ \gamma (t), (f \circ \gamma)' (t))= K_\Delta (\gamma (t), \gamma'(t))=1$$
we get
\begin{equation}\label{ppp}
(1-u)K_\Omega (f \circ \gamma_1(u), (f \circ \gamma_1)'(u))=1.
\end{equation}
There is a constant $C >0$ such that
\begin{equation}\label{kkk}
C^{-1} (1-u) \le \delta(\gamma_1(u)) \le C(1-u) \ \ \ \ \forall u \in [0,1].
\end{equation}
Also there exist $D>0$ with the following property:
\begin{equation}\label{lll}
D^{-1} \delta \le \delta (f \circ \gamma_1(u)) \le D \delta  \ \ \ \ \forall u \in [0,1].
\end{equation}
By Graham's estimates (see~\cite{gra}), there is an $E>0$ such that
$$ K_\Omega (p, v) \le E \frac {\Vert v \Vert}{\delta(p)}$$
for all $p \in \Omega$ and $v \in {\mathbb C}^n$.
Combining this estimate with (\ref{ppp}), (\ref{kkk}) and (\ref{lll}) we see that
$\lim_{u \rt 1}(f \circ \gamma_1)'(1) \neq 0$.

\hfill $\square$

 Let $U$ be a small neighborhood of $y^\infty$ in $\mathbb C^n$ and let $\Gamma \subset \Omega \cap U$ be a part of a half c\^one with vertex at $y^\infty$, axis tangent to $f \circ \gamma_1$ at $y^\infty$.

\begin{lemma}\label{bounded-lem}
There is a constant $C >0$ such that for every point $z \in f \circ \gamma([0,\infty)) \cap \Gamma$ we have :
$$
d^K_\Omega(z,\partial \Gamma) \geq C \sqrt{\delta(z)}.
$$
\end{lemma}

\noindent{\bf Proof :} 
Without loss of generality we may assume that $y^\infty = 0 \in \mathbb C^n$, that $z :(0',-\delta(z))$. Moreover there is a constant $c'>0$ such that $dist(z,\partial \Gamma) \geq c' \delta(z)$. Choosing $U$ sufficiently small we may assume that for every $p \in \Omega \cap U$ and for every $v \in \mathbb C^n$ (see \cite{gra})~:
$$
K_\Omega(p,v) \geq D\frac{\|v\|}{\sqrt{\delta(p)}}
$$
where $D > 0$.

In particular, consider a point $q \in \partial \Gamma$ such that $\delta(q) \leq 2 \delta(z)$ and a $C^1$ path $\alpha$ joining $z$ to $q$. We may restrict to the part of the path contained in $\Gamma \cap U$, implying that $\|\alpha(t)\| \sim \delta(\alpha(t))$, and in the ball centered at the origin, with radius $4 \delta(z)$.
Then there is a constant $c''>0$ such that :
$$
\displaystyle l(\alpha)  =  \int_0^1 K_\Omega(\alpha(t),\alpha'(t))dt  \geq  \displaystyle  c'' \int_0^1\frac{\|\alpha'(t)\|}{\sqrt{\|\alpha(t)\|}}dt.
$$

By our restriction we have :
$$
l(\alpha) \geq \frac{c''}{2\sqrt{\delta(z)}}\int_0^1\|\alpha'(t)\|dt \geq \frac{c''}{2\sqrt{\delta(z)}} \Vert z - q \Vert \geq D' \sqrt{\delta(z)}
$$
by definition of $\Gamma$. This proves Lemma \ref{bounded-lem}. The claim is now a direct consequence of Lemma \ref{bounded-lem} and of Condition (\ref{equiv-eq}). \vspace{2mm}

Finally to complete the proof of Lemma~\ref{lem-lambda}, we show that $\lambda=1$ where $\lambda$ is given by the equality
$(f \circ \gamma_1)'(1) = \lambda (f \circ \sigma_1)'(1)$.
\begin{lemma}\label{nbn}
Let $\Omega$ be a $C^2$ strongly pseudoconvex domain in $\C^n$. The function defined on $\Omega \times \C^n$ by
$$(z, v) \mapsto \delta(z) K_\Omega(z,v)$$
extends continuously to $\overline \Omega \times \C^n$.
\end{lemma}

\noindent{\bf Proof:} According to \cite{gra} we have $\lim_{z \to p \in \partial \Omega} \delta(z) K_\Omega(z,v) = \frac{1}{2}\|v_N(p)\|$ where $v_N(p)$ is the complex normal component to $\partial \Omega$ of the vector $v$ at $p$. The result follows since the map $(p,v) \mapsto v_N(p)$ depends continuously on $(p,v) \in \overline \Omega \times \mathbb C^n$. \qed
\vspace{2mm}

Lemma \ref{nbn} implies that
$$\lim_{u \rt 1} \delta ( f \circ \gamma_1 (u)) K_\Omega( f \circ \gamma_1(u),(f \circ \gamma_1)'(u))= \lambda
\lim_{u \rt 1} \delta ( f\circ \sigma_1 (u)) K_\Omega( f \circ \sigma_1(u),(f \circ \sigma_1)'(u)). $$


This can be written as
$$ \lim_{u \rt 1} \frac {\delta ( f \circ \gamma_1 (u)) K_\Omega( f \circ \gamma_1(u),(f \circ \gamma_1)'(u))} {\delta ( f\circ \sigma_1 (u))K_\Omega( f \circ \sigma_1(u),(f \circ \sigma_1)'(u))} = \lambda.$$

By (\ref{pop}) and the assumption that $\gamma$ and $\sigma$ have unit-speed, the above equality gives
$$ \lim_{t \rt \infty} \frac {\delta ( f \circ \gamma (t))}{\delta ( f\circ \sigma (t))} = \lambda.$$
We claim that the left-hand side above is equal to $1$. This follows from the following lemma.

\begin{lemma}\label{last-lem}
Let $\Omega$ be a $C^2$ strongly pseudoconvex domain in $\C^n$. Let $\{x_k\}, \ \{y_k\}$ be sequences in
$\Omega$ satisfying $\lim_{k \rt \infty}x_k = \lim_{k \rt \infty}y_k =z \in \partial \Omega$
and $\lim_{k \rt \infty} d_\Omega^K (x_k,y_k) =0$. Then
$$\lim_{k \rt \infty} \frac {\delta(x_k)}{\delta(y_k)}=1.$$

\end{lemma}

{\bf Proof:} This follows from the calculations in the proof of Lemma 4.1 of ~\cite{bb}. Let $\Omega$
be as above. Let $\pi: T \rt \partial \Omega$ be the closest point projection map, which is well-defined
in an $\epsilon$-tubular neighbourhood $T$ of $\partial \Omega$ for $\epsilon$ small enough. For $q \in \partial \Omega$, let
$n(q)$ denote the outer unit normal to $\partial \Omega$ at $q$ and for $p \in T \cap \Omega$ let
$$ h(p) =\sqrt {\delta(p)}.$$
If $\gamma:[0,1] \rt T \cap \Omega$ is a $C^1$ curve with $\gamma(0)=x_k$ and
$\gamma(1)=y_k$, then
\begin{align} \notag
\Bigl \vert \frac {d}{dt} h(\gamma(t)) \Bigr \vert &=  \frac {1}{2 h(\gamma(t))} \Bigl \vert \frac {d}{dt} \delta (\gamma(t)) \Bigr \vert = \frac {\vert Re \langle \bar \partial \delta(\gamma(t)), \gamma '(t) \rangle \vert}{h(\gamma(t))} \notag \\
&=   \frac {\vert Re \langle n (\pi (\gamma(t))), \gamma '(t) \rangle \vert}{2h(\gamma(t))} \le \frac {\vert \gamma_N '(t) \vert}{2 h(\gamma(t))}, \label{ggg}
\end{align}
where $\gamma _N ^{'}(t)$ is the normal component of $\gamma (t)$ in the
standard decomposition. There is a constant $C>0$ (depending only on $T$)  such that the Kobayashi length $l(\gamma)$ of $\gamma$ satisfies
\begin{align}\notag
l(\gamma) & \ge C \int_0^1  \frac {\vert \gamma_N '(t) \vert}{ h(\gamma(t))^2}dt \notag \\
& \ge 2C \int_0^1 (log(h(\gamma(t))))' dt =
2C \Bigl \vert log \Bigl ( \frac{h(y_k)}{h(x_k)} \Bigr ) \Bigr \vert \label{cdd}
\end{align}

Next let $\gamma$ be a curve in $\Omega$ connecting $x_k$ and $y_k$ which exits $T \cap \Omega$. Let
$\gamma(t_1), \gamma(t_2) \in \partial T \cap \Omega$ be the first exit point and last entry point respectively. We have $h(\gamma(t_1)) = h(\gamma(t_2))=\sqrt \epsilon$.  By the above estimate
$$ l( \gamma \vert_{[0,t_1]}) \ge 2C \Bigl \vert log \Bigl ( \frac{\sqrt \epsilon}{h(x_k)} \Bigr ) \Bigr \vert
\ \ \ \ \ {\rm and} \ \ \ \ \ l( \gamma \vert_{[t_2,1]}) \ge 2C \Bigl \vert log \Bigl ( \frac {h(y_k)}{\sqrt \epsilon} \Bigr ) \Bigr \vert .$$
Adding the two inequalities above we get
$$ l(\gamma) \ge l( \gamma \vert_{[0,t_1]}) + l( \gamma \vert_{[t_2,1]}) \ge 2C \Bigl \vert log \Bigl ( \frac{h(y_k)}{h(x_k)} \Bigr ) \Bigr \vert. $$

Hence (\ref{cdd}) holds for all $C^1$ curves connecting $x_k$ and $y_k$ and  we get
$$d_\Omega^K(x_k,y_k) \ge 2C \Bigl \vert log \Bigl ( \frac{h(y_k)}{h(x_k)} \Bigr ) \Bigr \vert .$$
This gives the desired conclusion. \hfill $\square$

\section{Proof of Theorem \ref{mai1} and of Theorem \ref{mai2}}

In order to prove Theorem \ref{mai1} we need to establish the following statement :

\begin{proposition}\label{discs-prop}
Let $f : \Delta \to \Omega$ be an isometry. Then $f$ is either holomorphic or anti-holomorphic.
\end{proposition}

\noindent{\bf Proof:}
Let $\gamma, \sigma : [0,\infty) \rightarrow \Delta$ be two geodesics as in Lemma \ref{lem-lambda}. In particular,
$\lim_{t \rightarrow \infty} \gamma(t) =\lim_{t \rightarrow \infty} \sigma (t)=1 \in \partial \Delta$.
Let $\gamma_1, \ \sigma_1: [0,1] \rt \Omega$ be the corresponding reparametrizations. Then
\begin{equation}\label{aqq}
(f \circ \gamma_1)'(1)=(f \circ \sigma_1)'(1)
\end{equation}
by Lemma \ref{lem-lambda}.  By Lemma \ref{geod-lem} there are complex geodesics $\phi_\gamma,\phi_\sigma : \Delta \to \Omega$ and real geodesics $\tilde \gamma$ and $\tilde \sigma$ in $\Delta$, with reparametrizations $\tilde \gamma_1$ and $\tilde \sigma_1$, such that
\begin{equation}\label{aqq2}
\phi_\gamma (\tilde \gamma_1(t)) = f (\gamma_1(t)), \ \ \phi_\sigma (\tilde \sigma_1(t)) = f (\sigma_1(t))
\end{equation}
for all $t \in [0,1]$. We may assume that $\tilde \gamma(1)=\tilde \sigma(1)=1 \in \partial \Delta$ and that $\tilde \gamma_1'(1)= \tilde \sigma_1'(1)=:v$. 
It follows from
(\ref{aqq}) and (\ref{aqq2}) that
$$ d \phi_\gamma(1)(v)= d \phi_\sigma(1)(v).$$
Both $\phi_\gamma$ and $\phi_\sigma$ being holomorphic it follows that $d \phi_\gamma(1) =d \phi_\sigma(1)$.


By Lemma~\ref{kol}, we have $\phi_\gamma = \phi_\sigma$. Fixing a geodesic $\gamma$ in $\Delta$, the set of all the geodesics $\sigma$ in $\Delta$ such that $\gamma$ and $\sigma$ satisfy the assumptions of Lemma \ref{lem-lambda} forms a foliation of $\Delta$. It follows in particular that $f(\Delta)=\phi_\gamma(\Delta)$.

Without loss of generality we can assume that the geodesic $\gamma$ maps to the real line  in $\Delta$ i.e. $\gamma: (-\infty,\infty) \to \mathbb R \cap \Delta$ with $p:=\gamma(0)=0$. We also consider the geodesic $i\gamma: (-\infty,\infty) \to i\mathbb R \cap \Delta$. Keeping the notations of Remark \ref{geod-rem} let $\mu$ be an automorphism of $\Delta$ (or the conjugate of an automorphism of $\Delta$) such that $f \circ \gamma =  \phi_\gamma \circ \mu \circ \gamma$. In particular $f(0)=\phi_\gamma \circ \mu(0)$. Let $\tilde \gamma$ be the unique geodesic such that $\phi_\gamma \circ \mu \circ \tilde \gamma = f \circ i\gamma$. Then $\phi_\gamma \circ \mu(\tilde \gamma (0)) = \phi_\gamma \circ \mu(0)$ and necessarily $\tilde \gamma(0)=0$ since $\phi_\gamma \circ \mu$ is an embedding. Hence $\tilde \gamma = e^{i\theta} \gamma$ for some $\theta \in \mathbb R$. 

Let $\lambda$ be the unique real geodesic in $\Delta$ joining $\gamma(1)$ and $i\gamma(1)$. The unique geodesic $\tilde \lambda$ such that $f \circ \lambda = \phi_\gamma \circ \mu \circ \tilde \lambda$ joins the two points $\gamma(1)$ and $\tilde \gamma(1)$. Hence we have $\lambda(0)=\tilde \lambda(0)=\gamma(1)$ and $\lambda(t_0)=i\gamma(1)$ for some $t_0 > 0$. Therefore $\tilde \lambda(t_0)=\tilde \gamma (1)$. We have now :
$$
\begin{array}{lllllllll}
d^K_\Delta(\gamma(1),i\gamma(1)) & = & d^K_\Delta(\lambda(0),\lambda(t_0)) \\
& & \\
& = & d^K_\Omega(f(\lambda (0)),f(\lambda(t_0))) \\
& & \\
& = & d^K_\Omega(\phi_\gamma \circ \mu(\tilde \lambda (0)),\phi_\gamma \circ \mu(\tilde \lambda(t_0))) \\
& & \\
& = & d^K_\Delta(\tilde \lambda(0),\tilde \lambda(t_0)) \\
& & \\
& = & d^K_\Delta(\gamma(1),\tilde \gamma(1))\\
& & \\
& = & d^K_\Delta(\gamma(1),e^{i\theta}\gamma(1)).
\end{array}
$$
It follows that $\theta \equiv \pi/2 \ mod [\pi]$. 

If $\theta \equiv \pi/2 \ mod [2\pi]$ (resp. $\theta \equiv 3 \pi/2 \ mod [2\pi]$) then  $f$ and $\phi \circ \mu$ (resp. $f$ and $ \phi \circ \overline \mu$) agree on the geodesics $\gamma$ and $i \gamma$.  If $z \in \Delta$ is any point then $z$ lies on a geodesic $\alpha: (-\infty, \infty) \rt \Omega$ passing through $\gamma(t_0)$ and $i \gamma(t_0)$ for some $t_0 \in (-\infty, \infty)$. By the uniqueness of real geodesics in $\Omega$ (Lemma \ref{geod-lem}), it follows that $f( \alpha(t))=\phi \circ \mu ( \alpha(t))$ (or $f \circ \overline \mu (\alpha(t))$) for all  $t \in (-\infty, \infty)$ since both are geodesics passing through $f (\gamma(t_0))$ and $f(i \gamma(t_0))$. Since $z$ lies on the image of $\alpha$, it follows that $f(z)=\phi \circ \mu (z)$ (or $f(z)= \phi \circ \overline \mu (z)$). \qed

\vskip 0,1cm
The proof of Theorem \ref{mai2} follows easily. Let $f: \Omega_1 \rt \Omega_2$ be a $C^1$ isometry. Let $p \in \Omega_1$
and let $J_i$ denote the almost-complex structure on $\R^{2n_i}, \ i=1,2$.  We want to prove that $J_2 \circ df_p =
df_p \circ J_1$. Let $V \subset T_p \Omega_1$ be a 2-dimensional  $J_1$ invariant subspace. We claim that $df_p(V)$ is $J_2$ invariant: Let $\phi: \Delta \rt \Omega_1$ be a complex geodesic with $T_p(\phi( \Delta)) = V$. Now $ f \circ \phi: \Delta \rt \Omega_2$ is again an isometry and hence holomorphic or anti-holomorphic by Theorem \ref{mai1}. Hence $df_p(V)$ is $J_2$ invariant. Since the restriction of the infinitesimal Kobayashi metric of a strongly convex domain  to a  2-dimensional $J$ invariant subspace is Hermitian (this follows from the existence of a complex geodesic tangent to the given subspace) it follows that if $V= Span_\R \{v,J_1(v)\}$ then $df_p(v)= \pm J_2(df_p(v))$. Hence $f$ is either holommorphic or anti-holomorphic at every point $p \in \Omega_1$. By continuity $f$ is either holomorphic everywhere or anti-holomorphic everywhere.




\vskip 0,5cm
{\small
\noindent Herv\'e Gaussier\\
(1) UJF-Grenoble 1, Institut Fourier, Grenoble, F-38402, France\\
(2) CNRS UMR5582, Institut Fourier, Grenoble, F-38041, France\\
{\sl E-mail address} : herve.gaussier@ujf-grenoble.fr\\
\\
Harish Seshadri\\
Department of Mathematics, Indian Institute of Science, Bangalore 560012, India\\
{\sl E-mail address} : harish@math.iisc.ernet.in
}

\end{document}